# A Method of Incorporating Matrix Theory to Create Mathematical Function-Based Music

**Sidarth Jayadev**

**Abstract**

This paper attempts to look for a mathematical method of composing music by incorporating Schonberg's idea of tone rows, and matrix theory from linear algebra. The elements of a "note-set" $S$ are considered as the integer values for the natural notes based on the C Major Scale and rational numbers for semitones. The elements of $S$ are effectively mapped by a polynomial function to another note-set $T$. To accomplish this, $S$ is treated as a column vector, applied to the matrix equation $A\mathbf{x} = \mathbf{b}$, where $\mathbf{x}$ denotes the vector $S$, $\mathbf{b}$ denotes the resulting set $T$, and A represents a square matrix. This method yields functions capable of mapping input note-sets to others, thereby creating collections of sets that can be permuted in any order to form musical harmonies.

**I. Historical Application of Sets in Composition**

The first recorded use of music set theory was by the celebrated composer Howard Hanson in 1960. He used the key concepts of music set theory to establish the connection between different musical "objects" such as pitches and rhythms (Hanson, 1960). Composers such as Arnold Schonberg and Allen Forte, however, largely established the connection between music set theory and atonal music. The main parts of music set theory originally dealt with permutations and combinations of pitches, formally called "pitch class sets" (Hanson, 1960). Pitch class sets were considered to be sets of notes, represented as integers that could be assembled in any way using techniques from a subject of mathematics called combinatorics (Forte, 1973).The pitch class set was either "well-tempered" (where the frequency ratio between each note was constant) or "non equal-tempered." Pitch class sets could be used to organize music in such a way that certain musical themes or motifs would be highlighted (Rahn, 1980).These sets were also used to allow certain transformations, such as inversion, transposition, and retrograde motion. This idea of using operations ultimately gave rise to the concept of pitch set "symmetry." This symmetry, in essence, basically resulted from the "identity operation" performed on the set itself (Rahn, 1980). This new idea of symmetry, identity operations, and transposition allowed composers to use the principles of set theory to create new and unique types of pieces – especially those categorized under 12-tone or atonal music pioneered by German composer, Arnold Schonberg. Although 12-tone theory was popular during Schonberg's time, it called for restrictions of notes and harmonies that could be used in a pitch set. A significant extension of Schonberg's 12-tone system was "serialism." Serial compositions incorporated elements from 12-tone theory, and also introduced new techniques for modifying pitch sets.

In this paper, a similar but unique method is proposed to apply set theory to music. Although the idea of a pitch set is incorporated (in the paper it is simply referred to as a "note-set"), in order to create compositions, the proposal is to use mathematical functions. A basic "note-set" is created, from which functions can map one set to another, generating new sets of notes that can be organized in any fashion. This application of set theory, unlike 12-tone theory, does not restrict the inclusion of certain notes and harmonies like the traditional pitch set. In Section II, a "note-set" will be considered, and mathematical functions are devised that map between different sets. In Section III, this new system will be used to compose musical examples where the mathematical foundation holds. Section IV, examines the potential of the theory, and its future applications to large pieces of music.

**II. Creation of Functions to map correspondences between note-sets**

The note-set will begin with the scale from $C_4$ to $C_5$, including all the natural notes and sharp/flat notes, depending on the enharmonic spelling. Each note will be represented by a number. $C_4 = 0$, $C^{\#}_4 = 1$, and so on, until $C_5 = 12$. If this convention is used, the note-set will be $\mathbf{S} = \{0, 1, 2, 3, 4, 5, 6, 7, 8, 9, 10, 11\}$. Coincidentally, this is also similar to the set used in the theory of serial music. The goal is to create a function, $f(n)$, that would map set $\mathbf{S}$ to another note-set $\mathbf{T}$, with different numerical and musical values. How can such a function be derived? What are the conditions under which the domain of this function would exist? In the system, notes that can be generated by discrete sounds are used; hence, integers can be substituted for these values. In the traditional serial system, set $\mathbf{S}$ contains exactly 12 distinct notes; the proposed system, on the other hand, is not limited to this restriction. Half-integers and quarter-integers are used to describe more notes such as semitones or quarter-tones. In fact, any real number can be used in general as long as the value is not irrational. Since an irrational number cannot be represented by a fraction, it does not have a definitive note value; therefore, values used in the note-set contain no irrational numbers. The domain of note-sets consists of the integers and rational numbers.

With this domain and range restriction, can a specific form of the function $f(n)$? Since the domain of *n* doesn't include the irrationals in this context system, one type of function that would satisfy this constraint is a polynomial of degree M with rational coefficients. Such a polynomial of degree M can be expressed in the following form:

$$f(n) = \sum_{i=0}^{M} c_i n^i$$

Now that there is a desired form of the function $f(n)$, the polynomial coefficients $c_i$ need to be determined. Given two note-set with the same number of elements, a one-to-one correspondence can be created between both sets. The first note-set is denoted by **S,** and the second note-set by **T**. The following can be stated, given that $s_i$ is an element of **S**, and $t_i$ is an element of **T**:

$$f(s_i) = t_i \qquad (1)$$

How can equation (1) help determine what $f(n)$ might be? Since each $n$ in $f(n)$ is $s_i$ and each $f(n)$ is $t_i$, a matrix (an array of numbers) of the possible input and output values can be created. Here, the following matrix equation will be used to deduce the coefficients of $f(n)$:

$$\begin{pmatrix} s_1^M & \cdots & 1 \\ \vdots & \ddots & \vdots \\ s_i^M & \cdots & 1 \end{pmatrix} \begin{pmatrix} c_i \\ c_0 \end{pmatrix} = \begin{pmatrix} t_1 \\ t_i \end{pmatrix} \qquad (2)$$

Equation (2) is commonly known in linear algebra as a "matrix equation." Simply stated, a matrix equation is an equation involving a square matrix (where the # of rows = # of columns) and two "column vectors" (columns with 1 row and *i* number of rows). The square matrix is conventionally called A, the first column vector **x**, and the second column vector **b**. The matrix equation reads A**x** = **b**. If the matrix equation can be solved, the correct coefficients for our polynomial function $f(n)$ can be calculated. If we find the coefficients, we can write a function that would correctly map between the first note-set **S** and the second set **T**. This function would then serve as the generalized formula for a polynomial of any degree that could map **S** to **T**. (This is one interpretation of a method known as Lagrange interpolation.) An example of a closed-form formula for the 4 x 4 matrix uses Cramer's Rule to determine all of the elements of the column vector **x**. This mathematical theorem can also be used to solve the general formula for a 12 x 12 matrix as well. An example of the first coefficient $c_3$ is shown below:

$$c_3 = \frac{b_4[n_2^2(n_3-n_4)-n_2(n_3^2-n_4^2)+(n_3^2n_4-n_3n_4^2)]-n_1^2[b_3(n_3-n_4)-n_2(b_2-b_1)+(b_2n_4-n_3b_1)]+n_1^1[b_3(n_3^2-n_4^2)-n_2^2(b_2-b_1)+(b_2n_4^2-n_3^2b_1)]-b_3(n_3^2n_4-n_3n_4^2)+n_2^2(b_2n_4-n_3b_1)-n_2(b_2n_4^2-n_3^2b_1)}{n_1^3[n_2^2(n_3-n_4)-n_2(n_3^2-n_4^2)+(n_3^2n_4-n_3n_4^2)]-n_1^2[n_2^3(n_3-n_4)-n_2(n_3^3-n_4^3)+(n_3^3n_4-n_3n_4^3)]+n_1^1[n_2^3(n_3^2-n_4^2)-n_2^2(n_3^3-n_4^3)+(n_3^3n_4^2-n_3^2n_4^3)]-n_2^3(n_3^2n_4-n_3n_4^2)+n_2^2(n_3^3n_4-n_3n_4^3)-n_2(n_3^3n_4^2-n_3^2n_4^3)}$$

**Note: The formula above is denoted by small font set to fit the denominator as to prevent the fraction from being decomposed into partial fractions.**

### III. Employing the Functions in Composition

A generalized method was developed in the previous section to find a specific polynomial mapping elements of note-sets **S** to **T**. In certain cases, this method can be used to map one note-set in a particular key to a note-set in a different key. In most cases, for 4 or 5-note-sets, this method will be applied to determine the function. It is important to note that depending the enharmonic spelling of certain notes, the *mathematical function* might differ from the corresponding *musical functions*. Suppose the note-set **S** is {0, 4, 7, 10}. Set **S** corresponds to {$C_4$, $E_4$, $G_4$, $B^ب_4$} or {$C_4$, $E_4$, $G_4$, $A^\#_4$} (depending on the enharmonic spelling). The first set would normally modulate to the "I chord" in F major. However, the second set, which represents an inversion of an Augmented $6^{th}$ chord, would progress to the "V chord" of E minor or E major. This difference would have an effect on the function, because the corresponding set would be different in each case, due to the different musical values for each resolution.

Consider a function mapping from one note-set to another. In these few examples the order of notes will be preserved as they are shown. The initial note-set, **S**, will be

{-10, -3, 0, 6}. Set **S** corresponds to the root position of the V chord in G Major: {$D_3$, $A_3$, $C_4$, $F^\#_4$}. If the linear function $f(n) = n + 4$ is applied to **S**, a new note-set, **T**, will be obtained, which is {-6, 1, 4, 10}. This corresponds to {$F^\#_3$, $C^\#_4$, $E_4$, $A^\#_4$}, the root position of the V chord in B Major/Minor. Since *linear* function was used to map set **S** to set **T**, the inversion of the second chord was preserved, and so were the intervals. This is expected from a linear function correspondence. Therefore, a conjecture is generalized: All linear functions of the form $f(n) = n + a$, where "a" is an arbitrary integer, will preserve the intervals, qualities, and inversions of any chord or sequence of notes described by a note-set.

In this example, instead of starting off with a known function, the function will be solved using the matrix method. The initial note-set will be **S** = {-7, -2, 2, 8} which corresponds to the 1st inversion of an Augmented 6th chord. The set **S** is *musically* equivalent to {$F_3$, $B^b_3$, $D_4$, $A^b_4$} and {$F_3$, $B^b_3$, $D_4$, $G^\#_4$}. The first note-set (with A-flat) resolves to the "I chord" in E-flat Major. In this scenario, the corresponding note-set is **T** = {$G_3$, $B^b_3$, $E^b_4$, $G_4$}. In numerical values, **T** = {-5, -2, 3, 7}. If a 4 x 4 matrix is used to solve for the coefficients of the desired polynomial function, we obtain the following:

$$f(n) = -\frac{47}{5400}n^3 + \frac{61}{5400}n^2 + \frac{3469}{2700}n + \frac{307}{645}.$$

This function would only work for mapping set **S** to set **T**, unless $f(n)$ returns other integers *not* found in sets **S** and **T**.

We will now consider that set **S** was spelled as {$F_3$, $B^b_3$, $D_4$, $G^\#_4$}. This would normally resolve to the "V chord" of D major/minor. In other words, set **T** = {$E_3$, $A_3$, $C^\#_4$, $A_4$} = {-8, -3, 1, 9}, according to standard counterpoint resolution. If the same matrix method is used to solve for the function, the following is obtained:

$$f(n) = \frac{1}{450}n^3 + \frac{7}{450}n^2 + \frac{223}{225}n - \frac{239}{225}.$$

This function is very different from the function previously solved for the "V – I" resolution. This demonstrates that depending on the resolution, a set with the same numerical values but different enharmonic note names often imply different functions for the corresponding set.

Consider a function algorithm, where a sequence of known functions is used to map from one note-set to another. Suppose set **S** has elements $\{C^\#_4, G_4, A_4, E_5, F^\#_5\} \rightarrow \{1, 7, 9, 16, 18\}$. The function $f(n) = n - 5$. Using $f(n)$, the set **T** is obtained - a perfect-5$^{th}$ transposition: **T** = {-4, 2, 4, 11, 13}. Consider a new function, $g(n) = -n + 6$. The function $g(n)$ maps **T** to **U** = {10, 4, 2, -5, -7}. If the function $h(n) = \frac{1}{2}n$ is applied, then set **V** = {5, 2, 1, -2.5, -3.5} is obtained. The musical equivalent is **V** = $\{E_4, D_4, C_4, A^{**}_4, A^*_4\}$. The note $A^{**}$, in this system, represents A ¾-sharp, while the note $A^*$ represents A ¼-sharp.

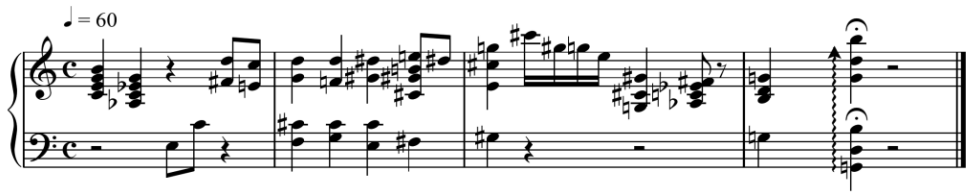

Figure 1. A short music composition sample demonstrating how mathematical functions can be used to create harmonies and melodies.

In this sample, three functions are involved, all of which are used in a row – similar to the function algorithm discussed above – to create a musically diverse harmonic progression. It is important to note that in the system, if a set contains a certain number of notes in any chosen order, the notes in the first octave must be used *before* the same notes in different octaves can be used. This is essential because in Figure 1, Bar 2, on the first beat, the note-set $\{F_3, C^\#_4, G_4, D_5\} \rightarrow \{-7, 1, 7, 14\}$ is arranged in its order. In the following beat, however, the order is switched – the G is an octave lower and F an octave higher – this switching is permitted *only after* the initial notes in their octave is used. The first function used to map the set {0, 4, 7, 11} (on beat 1, m.1) to {-4, 0, 3, 7} (on beat 2, m.1) is $f(n) = n - 4$. The corresponding note-sets for the first two beats of measure 1 are the **T** = $\{C_4, E_4, G_4, B_4\}$ and **U** = $\{A^b_3, C_4, E^b_4, G_4\}$. The second function used to map
{-4, 0, 3, 7} $\rightarrow \{A^b_3, C_4, E^b_4, G_4\}$ to {-8, 0, 6, 14} $\rightarrow \{E_3, C_4, F^\#_4, D_5\}$ is $g(n) = 2n$. In this scenario (beat 3, m.1), the set {-8, 0, 6, 14} is split among the first 3 eighth notes, with E and C in the bass on the first two 8$^{th}$ notes, and F$^\#$ and D as a chord on the third 8$^{th}$ note. This is allowed since the notes can be arranged in any order and rhythm. The third function mapping the set {-8, 0, 6, 14} to the set {-7, 1, 7, 5}→{$F_3$, $C^{\#}_4$, $G_4$, $D^{\#}_5$} is $h(n) = n + 1$. This mapping occurs on the transition between the 4$^{th}$ beat of measure 1 and the 1$^{st}$ beat of measure 2. The second beat of measure 2 features another permutation of the note-set {-7, 1, 7, 5}. Between the second and third beats of measure 2, the function used to map the set {-7, 1, 7, 5} to the note-set {-8, 1, 8, 15}→{$E_3$, $C^{\#}_4$, $G^{\#}_4$, $D^{\#}_5$} was

$$k(n) = \frac{-1}{924}n^3 + \frac{5}{1232}n^2 + \frac{1105}{924}n - \frac{35}{176}.$$

A notable difference between this function and the previous ones is that this function contains complicated rational coefficients, whereas the others used simple integer coefficients. The difference can be explained by the fact that one of the elements of the note-set {$F_3$, $C^{\#}_4$, $G_4$, $D^{\#}_5$} and {$E_3$, $C^{\#}_4$, $G^{\#}_4$, $D^{\#}_5$} is the same.

The next note-set after the third beat of measure 2 is **V** = {-7, 1, 7, 15, 11, 16}→ {$F^{\#}_3$, $C^{\#}_4$, $G^{\#}_4$, $D^{\#}_5$, $B_4$, $E_4$}. On the first beat of the third measure, a completely different note-set is introduced: the note-set {-4, 1, 4, 7}→{$G^{\#}_3$, $C^{\#}_4$, $E_4$, $G_4$}. The second and third beats of measure 3 represent a permutation of the elements of the set **A** = {-4, 1, 4, 7}. The function mapping this set to the next set on the fourth beat is $f(n) = n - 1$, from **A** to **B** = {-5, 0, 3, 6}→{$A^{b}_3$, $C_4$, $E^{b}_4$, $F^{\#}_4$}. This note is spelled as a $F^{\#}$ instead of a G-flat, since the next chord is the V chord of C major (G major). This harmonic progression represents a "German" Augmented 6$^{th}$ resolution.

The function used to complete the resolution is

$$f(n) = \frac{1}{99}n^3 + \frac{2}{99}n^2 + \frac{28}{33}n - 1.$$

The last chord is simply a reordering of the elements in the last set, {-6, -1, 2, 7}→ {$G_3$, $B_3$, $D_4$, $G_4$}. Figure 1 represents an algorithm consisting of all these functions used in a row, to compose a melody with harmony.

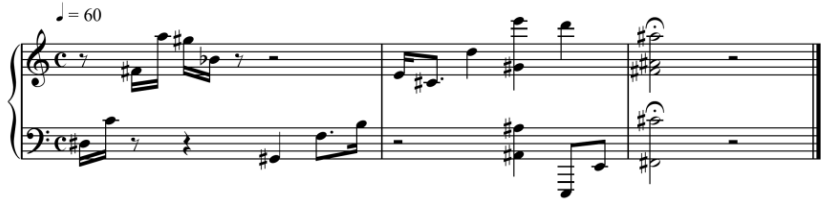

Figure 2. A shorter musical sample demonstrating function set mapping

In this sample (Figure 2), three functions were used to compose this sample. The first 2 beats represent the set $\mathbf{S} = \{21, 0, 6, -9, 11, 20\} \rightarrow \{A_5, C_4, F^\#_4, D^\#_3, B^b_4, G^\#_5\}$ prescribed in a certain order. The function $f(n) = n - 7$ was used to map $\mathbf{S}$ to the next note-set on the third and fourth beats, $\mathbf{T} = \{14, -7, -1, -16, 4, 13\} \rightarrow \{D_5, F_3, B_3, G^\#_2, E_4, C^\#_4\}$. The 3rd and 4th beats contain elements of the third note-set, $\mathbf{U} = \{28, -14, -2, -32, 8, 26\} \rightarrow \{E_6, A^\#_2, A^\#_3, E_2, G^\#_4, D_6\}$, mapped from $\mathbf{T}$ by the function $f(n) = 2n$. The second half of the fourth beat is the note $E_2$ transposed one octave higher, which is permitted in the system. The last function used to complete the melody mapped to the set $\{-6, 1, 6, 10, 17, 22\} \rightarrow \{F_3^\#, C^\#_4, A^\#_5, F^\#_4, F_5, A^\#_4\}$. Figure 2 – similar to Figure 1 - ends the continuous function algorithm used to create a melodic line with an harmonic progressions.

**IV. Special Chord Progressions and their associated function algorithms**

In practice, any arbitrary harmonic progression (represented by a chord progression from the root inversion) can be composed, and an algorithm can be associated as long as the cardinality of each note set is the same (so $f(x)$ is a bijection). In "traditional" harmony, a small set of progressions are often utilized repeatedly; here the function algorithms for mapping the essential chord qualities of each harmonic triad are categorized according to the key of each progression.

In C major, the basic (non-trivial) triad progression is the I → IV$^6_4$→V$^6$ → I (with an octave extension above the fifth of the triad). This corresponds to the following map:

$$\{0, 4, 7, 0\} \to \{0, 5, 8, 0\} \to \{11, 2, 7, 11\} \to \{0, 4, 7, 0\} \pmod{12}$$

By the matrix method for interpolating a polynomial across each pair-wise set, we obtain three polynomial functions, $f_1(x)$, $f_2(x)$, and $f_3(x)$, performed in that order, to generate this harmonic progression:

$$f_1(x) = -\frac{1}{28}x^2 + \frac{39}{28}x$$

$$f_2(x) = \frac{13}{30}x^2 - \frac{119}{30}x + 11$$

$$f_3(x) = -\frac{47}{180}x^3 + \frac{59}{20}x - \frac{77}{90}$$

For each of these functions, the highest order (degree) term is missing for a reason: specifically, it took on a free parameter, so to keep simplicity it was set equal to 0 for each function.

Here is another example of the same harmonic progression in D major:

$$g_1(x) = \frac{1}{84}x^2 + \frac{97}{84}x - \frac{5}{14}$$

$$g_2(x) = \frac{13}{180}x^2 - \frac{1}{20}x + \frac{73}{90}$$

$$g_3(x) = \frac{11}{120}x^2 + \frac{43}{24}x + \frac{3}{10}$$

Just as in the previous example, the free parameter for the highest-order variable was set to 0 for each function to keep simplicity of the functions. One of the first things to notice about each function (for both C and D major) is that the greatest common divisor of the denominators of each coefficient (for each individual function) is greater than 1, implying that none of the numbers are relatively prime. A less trivial observation is that there are at least 2 distinct greatest common divisors of all the denominator coefficients.

Another type of harmonic progression (involving three chords with octave extensions) is the "Neapolitan – V – I" progression. The actual progression is $II^{b}{}_{6} \rightarrow V^{6}{}_{4} \rightarrow I$.

In C major, this harmonic progression corresponds to

$$\{5, 8, 1, 5\} \to \{2, 7, 11, 2\} \to \{0, 4, 7, 0\} \pmod{12}$$

The corresponding polynomial functions in the algorithm are

$$h_1(x) = \frac{47}{84}x^2 - \frac{157}{28}x + \frac{337}{21}$$

$$h_2(x) = -\frac{1}{180}x^2 + \frac{17}{20}x - \frac{151}{90}$$

and for D Major, the 2 functions are

$$k_1(x) = -x^2 + 18x - 71$$

$$k_2(x) = \frac{59}{120}x^2 - \frac{121}{24x} + \frac{271}{20}$$

## V. Remarks on Mathematical Function-Based Composition

      Section III examined how functions can be used in an algorithm to create musical melodies, harmonies. The method of starting with an initial note-set **S**, applying polynomial functions in order to map to different note-sets, could potentially be used to create an entire composition. For example, with a certain theme and a note-set, an algorithm of polynomial functions can be used to write an entire symphony, or a concerto, by not only mapping to different note-sets, but also organizing the notes in any order as long as the reordering occurs *after* the initial presentation. This differs from the 20$^{th}$ century application of set theory to music, because the proposed system takes the creativity of music composition, and instead, allows the mathematics of functions to take over the act of writing the music. The mathematics, in effect, becomes the source of the music, while the composer rearranges the notes and devises rhythms to make his piece appalling to the ear. This function-based system

can be used not only to compose music from the perspective of tonal standards, but also atonal music–especially the 12-tone music style originally inspired by Schonberg.

Although this paper has addressed mathematical functions and their relationship to pitches, mathematics could also be applied in a similar manner to rhythmical patterns. Mathematical "operators" could be applied to functions in order to modify the rhythms of the notes in the initial note-set. The components of our system make it more accessible to use than Schonberg Theory – since these functions can also be used to imitate 12-tone music without adding the restriction of non-repetitive notes. The function-based system also shows how mathematics itself is connected to music.